\documentclass[11pt,a4paper,leqno]{article}
\usepackage{a4wide}
\usepackage{latexsym}
\usepackage{amsmath}
\usepackage{amssymb}
\usepackage{stackrel}  
\usepackage{color}
\usepackage{graphicx}
\usepackage[linesnumbered,ruled,vlined]{algorithm2e}

\usepackage{authblk}

\usepackage[hidelinks]{hyperref}

\newtheorem{defin}{Definition}

\newtheorem{prop}{Proposition}
\newtheorem{theo}{Theorem}
\newtheorem{corol}{Corollary}
\newtheorem{example}{Example}
\newenvironment{proof}{\medskip\par\noindent{\bf Proof}}{\hfill $\Box$
\medskip\par}

\newcommand{\sumS}[4]{\displaystyle\sum_{\nu={#1}}^\infty a_{#2}{#3}x^{#4}}
\newcommand{\sumSJ}[4]{\displaystyle\sum_{\nu={#1}}^\infty a_{#2}{#3}(x-1)^{#4}}

\newcommand{\lambdaQ}[2]{[#1]_q\left([#2]_q+2\right)}

\newcommand{\C}{\mathbb{C}}
\newcommand{\N}{\mathbb{N}}
\newcommand{\R}{\mathbb{R}}

\begin{document}
\title{On a moment generalization of some classical second-order differential equations generating classical orthogonal polynomials}

\author[1]{Edmundo J. Huertas}
\author[1]{Alberto Lastra}
\author[1]{V\'ictor Soto-Larrosa}
\affil[1]{Universidad de Alcal\'a, Dpto. F\'isica y Matem\'aticas, Alcal\'a de Henares, Madrid, Spain. {\tt edmundo.huertas@uah.es $\quad$ alberto.lastra@uah.es $\quad$ v.soto@uah.es}}


\date{}

\maketitle
\thispagestyle{empty}
{ \small \begin{center}
{\bf Abstract}
\end{center}

The aim of the work is to construct new polynomial systems, which are solutions to certain functional equations which generalize the second-order differential equations satisfied by the so called classical orthogonal polynomial families of Jacobi, Laguerre, Hermite and Bessel. These functional equations can be chosen to be of different type: fractional differential equations, q-difference equations, etc, which converge to their respective differential equations of the aforesaid classical orthogonal polynomials. In addition to this, there exists a confluence of both the families of polynomials constructed and the functional equations who approach to the classical families of polynomials and second-order differential equations, respectively
\smallskip

}
\bigskip \bigskip

\section{Introduction}\label{[SECTION-1]-Intro}

An infinite sequence of polynomials $(\varphi_n(x))_{n\ge0}$ with $\hbox{deg}(\varphi_n(x))=n$, is said to be orthogonal with respect to certain positive Borel measure $\mu$ supported in an infinite set $E\subseteq \R$, if all the following quantities, known as their \textit{moments}, satisfy 
$$\int_{E}x^nd\mu(x)<\infty,\quad n\in\N_0:=\{0,1,\ldots\},$$
and the polynomials in the sequence fulfill
$$\int_{E}\varphi_m(x)\varphi_n(x)d\mu(x)=\xi_n\delta_{m,n},\quad m,n\in\N_0,$$
where $\delta_{m,n}$ stands for the Kronecker delta (i.e. $\delta_{m,n}=0$ for $m\neq n$ and $\delta_{n,n}=1$), and where $\xi_n$ is a positive real number for all $n\in\N_0$ (see, for example \cite{Chi78, Ism05, Sze75}).

Said polynomial sequence $(\varphi_n(x))_{n\ge0}$ is also called \textit{classical}, if there exist a polynomial of degree at most 2, say $\sigma=\sigma(x)$ and a polynomial of degree at most one, say $\tau=\tau(x)$,  such that for every $n\in\N_0$ a real constant $\lambda_n$ exists for which $\varphi_n(x)$ satisfies the second-order differential equation 

\begin{equation}
    \sigma(x)y''(x)+\tau(x)y'(x) +\lambda_ny(x) = 0.
    \label{eq:SL}
\end{equation}

Observe that the restriction on the degree of the polynomials involved in (\ref{eq:SL}) guarantees that the solution to the Sturm-Liouville problem associated to the differential operator $\mathcal{L}[y](x) = \sigma(x)y''(x)+\tau(x)y'(x)$ has no irregular singularities. In 1929, S. Bochner~\cite{bochner} described all the families of equations of the form (\ref{eq:SL}) which admit a polynomial solution of degree $n$, for every $n\in\N_0$. He proved that, essentially, the only orthogonal polynomials satisfying this last property were those of Jacobi, Laguerre, Hermite and Bessel. Further classifications and details on orthogonal polynomials which solve second-order differential equations can be found in~\cite{kwonlittle} and the references therein.

In a different level, the concept of moment derivative was put forward by W. Balser and M. Yoshino in~\cite{bayo} in 2010, in the context of the study of formal solutions to functional equations. Moment differential equations (i.e. functional equations involving moment differential operator) have been proved to be of great versatility providing results which particularize to  classical differential equations, $q-$difference equations, or fractional differential equations when choosing adequately the sequence of moments. We refer to some of the recent advances in the study of moment differential equations described in~\cite{lamisu4,mi13,mi17,mitk,su} among others, regarding its original point of view, and~\cite{lastra,lastraprisuelos} in the context of the solution to systems of equations. 

The main aim in this work is to describe moment differential equations of second order, in the spirit of (\ref{eq:SL}), generalizing the classical ordinary differential equations associated to Jacobi, Laguerre, Hermite and Bessel polynomials, admitting a polynomial solution. Moreover, we will see that the classical equations and polynomials are recovered for a particular choice of the sequence of moments, whereas other choices provide families of second order functional equations of different nature (mainly fractional differential equations or $q$-difference equations) with polynomial-like or polynomial solutions. Such functional equations together with their solutions converge to the classical differential equations and the classical polynomials.

In a first section (Section~\ref{secmomder}), we have briefly described some preliminaries on moment derivation. Section~\ref{sec3} is devoted to the construction of certain families of polynomials generalizing Laguerre, Hermite, Jacobi and Bessel polynomials, satisfying a second order moment differential equation. The choice of two concrete moment sequence which are important in applications is made in Section~\ref{secapcon}, leading to polynomial solutions to fractional differential equations and $q$-difference equations. The work follows with some numerical results showing the confluence of the polynomials and functional equations to the classical ones. The work concludes with some conclusions and future directions of research.


\section{Preliminaries. Moment derivation}\label{secmomder}

In this preliminary section, we recall the main definitions on moment derivation, and moment differential operators with some of their central properties.

Let $m=(m(p))_{p\ge0}$ be a sequence of positive real numbers. The moment derivative operator $\partial_m$ is a formal operator defined on the vector space of formal power series with complex coefficients, $\C[[z]]$ given by
\begin{align*}
\partial_m:\C[[z]]&\to\C[[z]]\\
\sum_{p\ge0}\frac{a_p}{m(p)}z^p&\mapsto \sum_{p\ge0}\frac{a_{p+1}}{m(p)}z^p.
\end{align*}

It holds that 
$$\partial_m(x^p)=\left\{ \begin{array}{lc}
            0  & p=0 \\
            \frac{m(p)}{m(p-1)}x^{p-1} & p\ge1 
             \end{array}
   \right.$$

It is natural to extend the previous definition to holomorphic functions defined on some neighborhood of the origin by identification of the function with its Taylor expansion at 0. In~\cite{lamisu3,lamisu4,lamisu2}, the authors have also extended the definition of $\partial_m$ to functions which are generalized sums or generalized multisums of formal power series, in the sense of~\cite{lamasa,jikalasa}, respectively. We also refer to~\cite{sanz,sanzrev} for a broad sight on the theory regarding the development of kernel functions for generalized summability, as the fundamentals motivating the step forward from the classical theory (see Chapter 5,~\cite{ba2}). 

It is worth remarking that $m$ can be any sequence of positive real numbers. However, it is frequent that $m$ is the sequence of moments associated to some measure with support in $[0,\infty)$. For this reason, $\partial_m$ is usually known as the moment derivative. 

Some of the most important situations in applications are the following choices of the sequence $m$ that will be considered in detail in Section~\ref{secapcon}.

\section{Generalized moment analogs of classical equations and polynomials}\label{sec3}

In this section, we consider certain moment differential equations of second order generalizing the classical ODEs satisfied by the classical orthogonal polynomials. In the whole section, we fix a sequence of positive real numbers $m=(m(p))_{p\ge0}$, for which we assume it is normalized by $m(0)=1$. The main result of the present work provides generalizations of the classical second order differential equations, admitting polynomial or polynomial-like solutions, which satisfy a confluence to the classical equations and polynomials when the sequence $m$ approaches the sequence $(p!)_{p\ge0}$. 

\subsection{Generalized Laguerre polynomials}

Let $\alpha>-1$. For every integer $n\ge1$ we define $\lambda_{n}^{L}:=\frac{m(n)}{m(n-1)}$ and consider the moment differential equation
\begin{equation}
    x\partial^2_m y(x)+(\alpha+1-x)\partial_m y(x)+\lambda_n^{L}y(x) = 0
    \label{eq:MomentLaguerre}
\end{equation}

\begin{prop}\label{prop1}
There exists a polynomial solution $y_{L,m}$ of (\ref{eq:MomentLaguerre}) of degree $n$ given by
$$y_{L,m}(x)=\sum_{p=0}^{n}\frac{a_p}{m(p)}x^p$$
where 
$$a_1 = -\frac{\lambda_n^L}{\alpha+1}a_0,$$
and 
\begin{equation}
a_p = -\frac{m(p)\lambda_n^{L}}{m(1)(\alpha+1)}a_0\prod_{j=1}^{p-1}\left( \frac{\frac{m(j)}{m(j-1)}-\lambda_n^{L}}{\frac{m(j+1)}{m(j)}\left[\frac{m(j)}{m(j-1)}+\alpha+1\right]}\right)
\label{eq:RecurrenceMomentLaguerre}
\end{equation}
for $2\le p\le n$.
\end{prop}

\begin{proof}
Let us consider a formal power series $y(x)=\sum_{p\ge0}\frac{a_p}{m(p)}x^p$. By plugging the previous formal power series into (\ref{eq:MomentLaguerre}), one obtains the recursion formula for the coefficients $a_p$ after identification of the corresponding terms:
\begin{equation*}
    a_1 = \left[-\frac{\lambda_n^{L}}{\alpha+1}\right]a_0,
\end{equation*} 
and
\begin{equation*}
\frac{a_{p+1}}{m(p+1)}=\left[\frac{m(p)}{m(p+1)}\frac{\frac{m(p)}{m(p-1)}-\lambda_n^{L}}{\frac{m(p)}{m(p-1)}+\alpha+1}\right]\frac{a_p}{m(p)}
\end{equation*} 
for $p\geq 1$. Observe from the choice of $\lambda_n^{L}$ that $a_{j}=0$ for $j\ge n+1$, obtaining a polynomial of degree $n$. The other coefficients are determined by (\ref{eq:RecurrenceMomentLaguerre}), up to fixing $a_0$. 
\end{proof}

Observe that in the classical case, i.e. $m=(p!)_{p\ge0}$, the equation $(\ref{eq:SL})$ is given by
\begin{equation}
    xy''(x) + (\alpha+1-x)y'(x)+n y(x) = 0,
    \label{eq:Laguerre}
\end{equation}
which is recovered from (\ref{eq:MomentLaguerre}). Since $x = 0$ is a regular singular point of the equation, the well-known Frobenius method to search for formal solutions to the equation leads to a formal solution of (\ref{eq:Laguerre}) of the form 

\begin{equation}
y(x) = \sumS{0}{\nu}{}{\nu+r},
\label{eq:yLaguerre}
\end{equation}
for the values of $r$ satisfying the indicial equation $a_0r(r+\alpha)=0$, leading to the values $r_1 = 0$ and $r_2 = -\alpha$. The coefficients $a_\nu$ can be obtained by direct inspection when plugging the formal power series (\ref{eq:yLaguerre}) into (\ref{eq:Laguerre}). Indeed, one has that
\begin{equation*}
    a_\nu = \frac{(r-n)_\nu}{(r+1)_\nu(r+\alpha+1)_\nu}a_0, \quad \nu\geq 1
\end{equation*}
where $(x)_\mu$ stands for the Pochhammer symbol, and where $a_0\neq0$. Each of the choices for $r$ provides a solution. A polynomial solution of degree exactly $n$ is obtained for $r=0$. More precisely,
\begin{equation*}
y_1(x) = y(x)|_{r=0} = a_0\left(1+\displaystyle\sum_{\nu=1}^n  \frac{(-n)_\nu}{\nu!(\alpha+1)_\nu}x^\nu\right) = a_0\hspace{0.1cm}{_1F}_1(-n;\alpha+1;x).
\end{equation*}
Here, ${_1F}_1(a;b;x)$ stands for the confluent hypergeometric function.

Observe that $y_1$ turns out to be proportional to the generalized Laguerre polynomial $L^{(\alpha)}_n(x)$. 

\subsection{Generalized Hermite polynomials}

For every integer $n\ge1$ we define $\lambda_{n}^{H}:=\frac{2m(n)}{m(n-1)}$ and consider the moment differential equation
\begin{equation}
\partial_{m}^2y(x) -2x\partial_my(x)+\lambda_n^{H}y(x) =0.
\label{eq:MomentHermite}
\end{equation}
\begin{prop}\label{prop2}
There exists a polynomial solution $y_{H,m}$ of (\ref{eq:MomentHermite}) of degree $n$ given by
$$y_{H,m}(x)=\sum_{p=0}^{n}\frac{a_p}{m(p)}x^p,$$
where 
\begin{enumerate}
\item If $n$ is an even number, $a_{2k-1}=0$ for $1\le k\le n/2$ and
$$a_2=-\lambda^{H}_na_0,$$
\begin{equation}
a_{2k} =-m(2k)\frac{\lambda_n^{H}a_0}{m(2)}\prod_{j=1}^{k-1}\left(\frac{m(2j)\left[\frac{2m(2j)}{m(2j-1)}-\lambda_n^{H}\right]}{m(2j+2)}\right)
\label{eq:RecurrenceevenMomentHermite}
\end{equation}
for $k\ge2$.
\item If $n$ is an odd number, $a_{2k}=0$ for $0\le k\le  \left\lfloor n/2 \right\rfloor$ and
\begin{equation}
a_{2k+1}= m(2k+1)\frac{a_1}{m(1)}\prod_{j=0}^{k-1}\left(\frac{m(2j+1)\left[\frac{2m(2j+1)}{m(2j)}-\lambda_n^{H}\right]}{m(2j+3)}\right)
\label{eq:RecurrenceoddMomentHermite}
\end{equation}
for $k\ge1$.
\end{enumerate}
\end{prop}
\begin{proof}

Inserting a formal power series $y(x)=\sum_{p\ge0}\frac{a_p}{m(p)}x^p$ into equation (\ref{eq:MomentHermite}), one arrives at the recursion formula starting from $a_2 = -\lambda_n^{H}a_0$ for the even coefficients, and with
\begin{equation*}
\frac{a_ {p+2}}{m(p+2)} = \frac{m(p)\left[\frac{2m(p)}{m(p-1)} - \lambda_n^{H}\right]}{m(p+2)}\frac{a_{p}}{m(p)},
\end{equation*}
for $p\geq1$ in both cases. The result follows after collecting all the terms in the recurrence.

Observe from (\ref{eq:RecurrenceevenMomentHermite}), (\ref{eq:RecurrenceoddMomentHermite}), and the choice of $\lambda_n^{H}$ that $y_{H,m}$ turns out to be a polynomial of degree $n$, independently of $n$ being an odd or even number. 
\end{proof}

Regarding the classical case, i.e. $m=(p!)_{p\ge0}$, an analogous reasoning can be followed for the second-order differential equation,
\begin{equation}
y''(x) -2xy'(x) + 2ny(x) = 0
\label{eq:Hermite}
\end{equation}
which has an irregular singularity at $\infty$. One can search for solutions of the form
\begin{equation*}
 y(x) = \sumS{0}{\nu}{}{\nu},
 \label{eq:yHermite}
\end{equation*}
leading to the recurrence 
\begin{equation}
a_{\nu+2} = \frac{2(\nu-n)}{(\nu+2)(\nu+1)} a_\nu
\label{eq:anHer}
\end{equation}
for $\nu\geq0$. We observe that fixing $a_1=0$, the recurrence determines that the series will contain only even powers of $x$ or, similarly, for $a_0 = 0$ it will contain only odd powers of $x$. The general solution of (\ref{eq:Hermite}) is a linear combination of the even and odd solutions. Indeed, the solution consisting of even powers is determined by
\begin{equation*}
y_{even}(x) = \displaystyle\sum_{k=0}^{\infty} a_{2k}x^{2k} = a_0\left(1+\displaystyle\sum_{k=1}^{n}\frac{2^{2k}(-n)_k}{(2k)!}x^{2k}\right)  = a_0\hspace{0.1cm} {_1F}_1\left(-n;\frac{1}{2};x^2\right),
\end{equation*}
whereas the solution consisting of odd powers is given by
\begin{equation*}
y_{odd}(x) = \displaystyle\sum_{k=0}^\infty a_{2k+1}x^{2k+1} = a_1\left(x+ \displaystyle\sum_{k=1}^n \frac{2^{2k}\left(-n\right)_k}{(2k+1)!}x^{2k+1}\right)=a_1x\hspace{0.1cm}{_1F}_1\left(-n;\frac{3}{2};x^2\right)
\end{equation*}
The coefficients of both solutions are obtained by taking the recurrence (\ref{eq:anHer}) back to the first term.

\vspace{0.3cm}

\textbf{Remark:} Observe that the parity of the generalized Hermite polynomials coincides with that of the classical Hermite polynomials.

\subsection{Generalized Jacobi polynomials}

Let $\alpha,\beta>-1$. For every integer $n\ge2$ we define 
\begin{equation}\label{e367}
\lambda_n^{J} = \frac{m(n)}{m(n-1)}\left[\frac{m(n-1)}{m(n-2)}+\alpha+\beta+2\right]
\end{equation}
and consider the moment differential equation
\begin{equation}
    (1-x^2)\partial_{m}^2y(x) +\left[\beta-\alpha-(\alpha+\beta+2)x\right]\partial_{m}y(x)+\lambda_n^{J}y(x) = 0
    \label{eq:MomentJacobi}
\end{equation}

\begin{prop}\label{prop3}
There exists a polynomial solution $y_{J,m}$ of (\ref{eq:MomentJacobi}) of degree $n$ given by
$$y_{J,m}(x)=\sum_{p=0}^{n}\frac{a_p}{m(p)}x^p,$$
where the sequence $(a_p)_{p\ge0}$ satisfies the following recursion formula for given $a_0,a_1$:
\begin{equation*}
a_2 = -(\beta-\alpha)a_1-\lambda_n^{J} a_0,\qquad a_3 = [m(1)(\alpha+\beta+2)+(\beta-\alpha)^2-\lambda_n^{J}]a_1+\lambda_n^{J}(\beta-\alpha)a_0,
\end{equation*}
and
\begin{multline*}
\frac{a_{p+2}}{m(p+2)} = \left[\frac{m(p)}{m(p+2)}\left(\frac{m(p)}{m(p-1)}\left[\frac{m(p-1)}{m(p-2)}+\alpha+\beta+2\right]-\lambda_n^{J}\right)\right]\frac{a_p}{m(p)}\\
-\left[(\beta-\alpha)\frac{m(p+1)}{m(p+2)}\right]\frac{a_{p+1}}{m(p+1)}
\end{multline*}
for $p\ge2$.
\end{prop}
\begin{proof}
It is a direct consequence of the recursion formula obtained for the coefficients of a formal solution of the equation (\ref{eq:MomentJacobi}). Observe that setting $a_{n+1}=0$ ($p=n-1$ in the recurrence) one arrives for $p=n$ at 
\begin{equation*}
\frac{a_{n+2}}{m(n+2)} = \left[\frac{m(n)}{m(n+2)}\left(\frac{m(n)}{m(n-1)}\left[\frac{m(n-1)}{m(n-2)}+\alpha+\beta+2\right]-\lambda_n^{J}\right)\right]\frac{a_n}{m(n)},
\end{equation*}
which yields $a_{n+2}=0$ from the definition of $\lambda_n^{J}$ in (\ref{e367}). The three term recursion formula allows us to conclude that $a_p=0$ for $p\ge n+1$.
\end{proof}

We remark that in the classical case, i.e. $m=(p!)_{p\ge0}$ one arrives at the differential equation
\begin{equation}
(1-x^2)y''(x) +[\beta-\alpha-(\alpha+\beta+2)x]y'(x)+\gamma(n)y(x) =0
\label{eq:Jacobi}
\end{equation}
with $\gamma(n) = n(n+\alpha+\beta+1)$ and $\alpha,\beta > -1$.
This equation has regular singular points at $x=\pm1$ so it is natural to look for power series solutions around $x=1$ using the Frobenius method, say
\begin{equation}
y(x) = \sumSJ{0}{\nu}{}{\nu+r}
\label{eq:yJacobi}
\end{equation}

The substitution of the formal series (\ref{eq:yJacobi}) into (\ref{eq:Jacobi}) determines the roots 
\begin{equation}\label{e248}
r_1 = 0,\quad r_2 =-\alpha 
\end{equation}
of the indicial equation, together with the recurrence  
\begin{equation*}
a_\nu = \frac{\gamma(n)-(\nu+r-1)(\nu+r+\alpha+\beta)}{2(\nu+r)(\nu+r+\alpha)}a_{\nu-1}=\frac{(n-\nu-r+1)(n+\alpha+\beta+r+\nu)}{2(\nu+r)(\nu+r+\alpha)}a_{\nu-1},
\end{equation*}
for $\nu \geq 1$. Each of the values in (\ref{e248}) contributes with a solution to (\ref{eq:Jacobi}). A polynomial solution is determined by $r_1=0$

\begin{equation*}
    y_1(x) =  a_0\hspace{0.1cm}{_2F}_1\left(-n,n+1+\alpha+\beta;1+\alpha;\frac{1}{2}(1-x)\right)
\end{equation*}

where $_2F_1(a,b;c;x)$ is the Gaussian hypergeometric series. A power series solution around $x_0=0$ can also be performed. Following analogous steps we find that a formal solution of the form
\begin{equation*}
    y(x) = \sumS{0}{\nu}{}{\nu}
\end{equation*}
for (\ref{eq:Jacobi}) satisfies the recurrence formula
\begin{equation*}
a_{\nu+2} = \frac{\nu(\nu+\alpha+\beta+1)-\gamma(n)}{(\nu+2)(\nu+1)}a_\nu - \frac{(\beta-\alpha)}{\nu+2}a_{\nu+1},
\end{equation*}
for $\nu\ge0$. As for the previous cases, this recurrence is directly obtained by plugging the formal power series into the equation, determining a convergent power series due to $x=0$ is an ordinary point of the equation under study.

Observe that, in order that the formal solution defines a polynomial of degree $n$, then we set $a_{n+1}=0$, which correspond to $\nu = n-1$. Then, from the eigenvalue of the Jacobi differential equation $\gamma(n) = n(n+\alpha+\beta+1)$ it is clear that for $\nu = n$ we have $a_{n+2}=0$. 

\vspace{0.3cm}

\textbf{Remark:} Observe that from the nature of the equation (\ref{eq:Jacobi}), we have written the formal solution as a formal power series centered at $x=1$. This can also be done when dealing with formal solutions of moment differential equations without any further technical difficulty. First, computing the formal solution of the moment differential equation in powers of $x$ and then rewriting the formal power series in powers of $x-1$ taking into account that 
$$x^{n}=\sum_{j=0}^{n}\binom{n}{j}(x-1)^{j},$$ for every $n\ge0$. 

\vspace{0.3cm}

\textbf{Remark:} Observe that the choice $\alpha=\beta$ provides an odd (resp. even) polynomial whenever $n$ is odd (resp. even), in the same fashion as in the classical Jacobi polynomial.

\subsection{Bessel polynomials}

For every integer $n\ge1$ we define the moment differential equations
$$  (x+1)\partial_m y(x)-m(1) y(x) =0,$$
for $n=1$, and for $n\ge2$ we put  
$$\lambda^{B}_n = \frac{m(n)}{m(n-1)}\left[\frac{m(n-1)}{m(n-2)}+2\right]$$  
and consider the moment differential equation
\begin{equation}\label{eq:BesselMDE}
x^2\partial_m^2 y(x) + 2(x+1)\partial_my(x) -\lambda^{B}_n y(x)=0,
\end{equation}

\begin{prop}\label{prop4}
There exists a polynomial solution $y_{B,m}$ of (\ref{eq:BesselMDE}) of degree $n$  given by
$$y_{B,m}(x)=a_0\left(1+x\right),\quad \hbox{for }n=1,$$
$$y_{B,m}(x) = a_0\left(1+\frac{ (2+m(1)) m(2)}{2 m(1)^2}x+\frac{ (2+m(1)) (-2 m^2(1)+(2+m(1)) m(2))}{4 m^2(1)}x^2\right), \quad \hbox{for }n=2,$$
\begin{align*}
    y_{B,m}(x) =  & a_0\left[1+\frac{\lambda^{B}_n}{2m(1)}x+\frac{\lambda^{B}_n}{4m(2)}\left(\lambda^{B}_n-2m(1)\right)x^2\right.\\
    &\left.+\frac{\lambda^{B}_n}{4m(2)}\left(\lambda^{B}_n-2m(1)\right)\displaystyle\sum_{p=3}^{n}\left(\displaystyle\prod_{j=2}^{p-1}\frac{\lambda^{B}_n-\frac{m(j)}{m(j-1)}\left[\frac{m(j-1)}{m(j-2)}+2\right]}{\frac{2m(j+1)}{m(j)}}\right)x^p\right],
\end{align*}
	for $n\ge3$.
\end{prop}
\begin{proof}
The result is clear for $n=1$. Let $n\ge2$ and consider a formal solution of (\ref{eq:BesselMDE}) in the form $y(x)=\sum_{p\ge0}\frac{a_p}{m(p)}x^p$. The following recurrence for the coefficients is obtained:
$$ a_1 = \frac{\lambda^{B}_n}{2}a_0,\quad a_2 = \frac{\lambda^{B}_n}{4}\left(\lambda^{B}_n-2m(1)\right)a_0,$$
and 
\begin{equation*}
    \frac{a_{p+1}}{m(p+1)} = \left[\frac{\lambda^{B}_n-\frac{m(p)}{m(p-1)}\left[\frac{m(p-1)}{m(p-2)}+2\right]}{\frac{2m(p+1)}{m(p)}}\right]\frac{a_p}{m(p)}
\end{equation*}
for $p\ge2$. We observe that the choice of $\lambda^{B}_n$ implies that $a_{n+1} = 0$, obtaining a polynomial solution of at most degree $n$. The general form of any coefficient in the polynomial solution is obtained by iteration of the recursion defining its terms.
\end{proof}

For the situation in which $m=(p!)_{p\ge0}$ one recovers the classical equation satisfied by Bessel polynomials,
$$  x^2y''(x)+2(x+1)y'(x)-n(n+1)=0,$$
In this case we have a double regular singular point at $x=0$, so applying again the Frobenius method we obtain a hypergeometric solution of the above differential equation for $r=0$ given by 
\begin{equation*}
    y(x) =  a_0\hspace{0.1cm}{_2}F_0\left(-n,n+1;-\frac{x}{2}\right)
\end{equation*}
which defines the Bessel polynomials for a suitable normalization.

\section{Applications and numerical confluence results}\label{secapcon}

In this section, we describe two particular situations of the previous section to certain choices of $m$ which are important in applications, namely fractional differential equations and $q$-difference equations, together with numerical confluence results for both of them, polynomials and equations, to the classical ones.

Both applications are derived from the main result of the present work.

\subsection{Fractional differential equations}
The Caputo fractional derivative of a positive order is defined as follows, see~\cite{olsp} as a classical reference.

\begin{defin}
Let $\mu>0$, $t>a$ be real numbers. The Caputo fractional derivative of order $\mu$ is
\begin{equation*}
{}^{C}_aD^\mu_tf(t):=\left\{
\begin{aligned}
&\frac{1}{\Gamma(n-\mu)}\displaystyle\int_{a}^t \frac{f^{(n)}(\tau)}{(t-\tau)^{\mu+1-n}}d\tau, \quad n-1<\mu<n\in\mathbb{N},\\
&\frac{d^n}{dt^n}f(t), \quad \mu = n \in \mathbb{N}, 
\end{aligned}
\right.
\end{equation*}
$\Gamma(\cdot)$ stands for Gamma function.
\end{defin}

In the following we consider $a=0$ and write ${}^{C}_aD_x^\mu \equiv {}^{C} \! D_x^\mu$ for simplicity. It is straight to check that
\begin{equation*}
{}^{C}D^\mu_t t^s:=\left\{
\begin{aligned}
&\frac{\Gamma(s+1)}{\Gamma(s-\mu+1)}t^{s-\mu}, \quad n-1<\mu<n\in\mathbb{N}, p>n+1, p\in \mathbb{R}\\
&0, \quad  n-1<\mu<n\in\mathbb{N}, p\leq n+1, p\in \mathbb{N}, 
\end{aligned}
\right.
\end{equation*}

Let $k>0$ be an integer, and consider the sequence $m_{1/k}:=(\Gamma(1+\frac{p}{k}))_{p\ge0}$. We observe that 
$$\Gamma\left(1+\frac{p}{k}\right)=k\int_0^{\infty}x^{p}x^ke^{-x^{k}}dx,$$
thus $m_{1/k}$ is a sequence of moments. In view of the previous properties, the formal differential operator $\partial_{m_{1/k}}$ satisfies
$$(\partial_{m_{1/k}}f)(x^{1/k})={}^{C} \! D_x^{1/k}(f(x^{1/k})),\qquad f\in\C[[x]].$$

In particular, observe that with the choice $k=1$, the classical derivative for $m_1=(p!)_{p\ge0}$ is recovered. Taking this last property into account, one arrives at the following realization of the results derived in Section~\ref{sec3}.

\begin{corol}
Let $\mu\in\mathbb{Q}_+$ and $\alpha>-1$. For every positive integer $n$ there exists a fractional differential equation
\begin{equation*}
x^\mu \left({}^{C} \! D_x^\mu\right)^2 y(x)+(\alpha+1-x^\mu){}^{C} \! D_x^\mu y(x)+ \frac{\Gamma(\mu)}{B(1+\mu(n-1),\mu)} y(x) = 0
\label{eq:CaputoLaguerre}
\end{equation*}
with solution given by the polynomial in $x^{\mu}$ of degree $n$ defined by
$$y_{L,m_{\mu}}(x)=\sum_{\nu=0}^{n}a_{\nu}x^{\mu\nu},$$
with $a_1=-\frac{a_0}{(\alpha+1)\mu B(1+\mu(n-1),\mu)} $
and
$$a_\nu =-\frac{a_0}{(\alpha+1)\mu B(1+\mu(n-1),\mu)}\prod_{j=1}^{\nu-1}\left(  \frac{B(1+\mu  (n-1),\mu )-B(1+\mu  (j-1),\mu )}{\frac{B(1+\mu  (n-1),\mu ) B(1+\mu  (j-1),\mu ) }{B(1+\mu  j,\mu )}\left(\alpha
   +1+\frac{\Gamma (\mu )}{B(1+\mu  (j-1),\mu )}\right)}\right)$$
for all $\nu\ge 2$. Here, $B(\cdot,\cdot)$ stands for the beta function.
\end{corol}

\begin{corol}
Let $\mu\in\mathbb{Q}_+$. For every positive integer $n$ there exists a fractional differential equation
\begin{equation*}
    \left({}^{C} \! D_x^\mu\right)^2 y(x) - 2x^\mu \left({}^{C} \! D_x^\mu\right) y(x) + \frac{2\Gamma(\mu)}{B(1+\mu(n-1),\mu)} y(x) =0
    \label{eq:CaputoHermite}
\end{equation*}
with polynomial solutions in $x^{\mu}$ of degree $n$, defined by
$$y_{H,m_{\mu}}(x)=\sum_{\nu=0}^{n}a_{\nu}x^{\mu\nu},$$
where:
\begin{enumerate}
\item If $n$ is an even number, then $a_{2k-1}=0$ for $1\le k\le n/2$,
$$a_2=-\frac{2\Gamma(\mu)}{B(1+\mu(n-1),\mu)\Gamma(1+2\mu)}a_0,$$
and
\begin{align*}
a_{2k}= -&\frac{2\Gamma(\mu)a_0}{B(1+\mu(n-1),\mu)\Gamma(1+2\mu)} \\ \times
&\prod_{j=1}^{k-1} \frac{2\Gamma(\mu)B(1+2j\mu,2\mu)\left[B(1+\mu(n-1),\mu)-B(1+\mu(2j-1),\mu)\right]}{\Gamma(2\mu)B(1+\mu(n-1),\mu)B(1+\mu(2j-1),\mu)}
\label{eq:anevenCaputoHermite}
\end{align*}
for $k\ge2$.
\item If $n$ is an odd number, then $a_{2k}=0$ for $0\le k\le \left\lfloor n/2 \right\rfloor$, and
$$a_{2k+1} = a_1\prod_{j=0}^{k-1}\frac{2\Gamma(\mu)B(1+(2j+1)\mu,2\mu)\left[B(1+\mu(n-1),\mu)-B(1+2j\mu,\mu)\right]}{\Gamma(2\mu)B(1+\mu(n-1),\mu)B(1+2j\mu,\mu)},$$
for $k\ge1$.
\end{enumerate}
\end{corol} 

\begin{corol}
Let $\mu\in\mathbb{Q}_+$ and let $\alpha,\beta>-1$. For every integer $n\ge2$ there exists a fractional differential equation
\begin{multline*}
(1-x^{2\mu})({}^{C}\! D_x^\mu)^2 y(x) \\
+[\beta-\alpha -(\alpha+\beta+2)x]{}^{C} \! D_x^\mu y(x) + \frac{\Gamma(\mu)\left[\Gamma(\mu)+(\alpha+\beta+2)B(1+\mu(n-2),\mu)\right]}{B(1+\mu(n-1),\mu)B(1+\mu(n-2),\mu)} y(x) = 0
\end{multline*}
with a solution given by the polynomial in $x^{\mu}$ of degree $n$ defined by
$$y_{J,m_{\mu}}(x)=\sum_{\nu=0}^{n}a_{\nu}x^{\mu\nu},$$
with
\begin{equation*}
a_2 = -\frac{(\beta-\alpha)\Gamma(1+\mu)}{\Gamma(1+2\mu)}a_1 -\frac{\Gamma(\mu)\left[\Gamma(\mu)+(\alpha+\beta+2)B(1+\mu(n-2),\mu)\right]}{\Gamma(1+2\mu)B(1+\mu(n-1),\mu)B(1+\mu(n-2),\mu)} a_0
\end{equation*}
\begin{multline*}
a_3=\frac{\Gamma(1+\mu)\left(\Gamma(1+\mu)(\alpha+\beta+2)+(\beta-\alpha)^2-\frac{\Gamma(\mu)\left[\Gamma(\mu)+(\alpha+\beta+2)B(1+\mu(n-2),\mu)\right]}{B(1+\mu(n-1),\mu)B(1+\mu(n-2),\mu)}\right)}{\Gamma(1+3\mu)}a_1\\
+\frac{\Gamma(\mu)(\beta-\alpha)\left[\Gamma(\mu)+(\alpha+\beta+2)B(1+\mu(n-2),\mu)\right]}{\Gamma(1+3\mu)B(1+\mu(n-1),\mu)B(1+\mu(n-2),\mu)}a_0
\end{multline*}
and
\begin{multline*}
a_{\nu+2}=\resizebox{0.99\textwidth}{!}{$\left[\frac{B(1+\mu\nu,\mu)B(1+\mu(\nu+1),\mu)}{\Gamma(\mu)}\left( \frac{\Gamma(\mu)+(\alpha+\beta+2)B(1+\mu(\nu-2),\mu)}{B(1+\mu(\nu-1),\mu)B(1+\mu(\nu-2),\mu)}- \frac{\Gamma(\mu)+(\alpha+\beta+2)B(1+\mu(n-2),\mu)}{B(1+\mu(n-1),\mu)B(1+\mu(n-2),\mu)}\right)\right]$}a_{\nu}\\\\
-\frac{B(1+\mu(\nu+1),\mu)(\beta-\alpha)}{\Gamma(\mu)}a_{\nu+1}
\end{multline*}
for $\nu\ge 2$.
\end{corol}

\begin{corol}
Let $\mu\in\mathbb{Q}_+$. For every integer $n\ge1$ there exists a fractional differential equation 
$$(x^{\mu}+1){^C}D_x^{\mu}y(x)-\Gamma(1+\mu)y(x)=0,$$
for $n=1$ and
\begin{equation*}
    x^2\left({^C} \! D_x^{\mu}\right)^2y(x) +2(x^{\mu}+1)\left({^C} \! D_x^{\mu}\right)y(x) +\frac{\Gamma(\mu)\left[\Gamma(\mu)+2B(1+(n-2)\mu,\mu)\right]}{B(1+(n-1)\mu,\mu)B(1+\mu(n-2),\mu)}y(x)=0,
    \label{eq:F-Bessel}
\end{equation*}
for $n\ge2$, with solutions given by the polynomial in $x^{\mu}$ of degree $n$ defined by
$$y_{B,m_\mu}(x)=\sum_{\nu=0}^{n}a_{\nu}x^{\mu\nu},$$
where
$$a_1=\frac{a_0}{\mu}
\frac{\Gamma(\mu)+2B(1+\mu(n-2),\mu)}{B(1+\mu(n-1),\mu)B(1+\mu(n-2),\mu)}
$$
$$a_2=
\frac{\Gamma^2(\mu)\left(\Gamma(\mu)+2B(1+\mu(n-2),\mu)\right)\left(\Gamma(\mu)+2\mu B(1+\mu(n-2),\mu)\left[\frac{1}{\mu}-B(1+\mu(n-1),\mu)\right]\right)}{4\Gamma(1+2\mu)B^2(1+\mu(n-1),\mu)B^2(1+\mu(n-2),\mu)}a_0
$$
and
\begin{multline*}   
a_{\nu}=\frac{\Gamma^2(\mu)\left(\Gamma(\mu)+2B(1+\mu(n-2),\mu)\right)\left(\Gamma(\mu)+2\mu B(1+\mu(n-2),\mu)\left[\frac{1}{\mu}-B(1+\mu(n-1),\mu)\right]\right)}{4\Gamma(1+2\mu)B^2(1+\mu(n-1),\mu)B^2(1+\mu(n-2),\mu)}a_0\\
\times\left(\displaystyle\prod_{j=2}^{\nu-1}\frac{B(1+\mu j,\mu)}{2}\left[\frac{\Gamma(\mu)+2B(1+\mu n,\mu)}{B(1+\mu(n-1),\mu)B(1+\mu(n-2),\mu)}-\frac{\Gamma(\mu)+2B(1+\mu j,\mu)}{B(1+\mu(j-1),\mu)B(1+\mu(j-2),\mu)}\right]\right)
\end{multline*}
for $\nu\ge3$.
\end{corol}

\subsection{$q$-difference equations}\label{sec42}

Let $q>0$ with $q\neq 1$ and consider the sequence $m_q=([p]_q!)_{p\ge0}$, where $[p]_q!$ stands for the $q$-factorial defined by $[0]_q!=1$ and for every positive integer $p$,  $[p]_q!=:=[p]_{q}\cdot [p-1]_{q}\cdot\ldots\cdot [1]_{q}$, where $[j]_q$ is the $j$-th $q$-number $[j]_q=1+q+\ldots+q^{j-1}=\frac{1-q^{j}}{1-q}$.

It is straight to check that the $q$-derivative defined in 1908 by F. H. Jackson in~\cite{jackson} 
$$D_{q,x}f(x)=\frac{f(qx)-f(x)}{qx-x}$$
formally coincides with the moment derivation $\partial_{m_q}$. We observe that 
$$D_{q,x}x^n=\frac{1-q^n}{1-q}x^n=[n]_qx^{n-1}$$
for every positive integer $n$.

Assume that $q>1$. In this second realization of a moment derivation, it holds that the sequence of $q$-factorials $m_q$ is quite related to the sequence $(q^{p(p-1)/2})_{p\ge0}$ in the sense that both sequences generate the same functional space. This last sequence is indeed a sequence of moments associated to different measures supported in $[0,\infty)$, defined by functions such as $\sqrt{2\pi\ln(q)}\exp(\frac{\ln^2(\sqrt{q}x)}{2\ln(q)})$ or $\ln(q)\Theta_{1/q}(x),$ where $\Theta_{1/q}$ stands for Jacobi Theta function 
$$\Theta_{1/q}(z)=\sum_{p\in\mathbb{Z}}q^{-\frac{p(p-1)}{2}}z^p,$$
which is a holomorphic function in $\C\setminus\{0\}$ with an essential singularity at the origin.   

In this context, the results of Section~\ref{sec3} read as follows.

\begin{corol}
Let $q\in\R$ with $q>0$ and $q\neq1$. We also fix $\alpha>-1$. For every integer $n$ there exists a $q$-difference equation of the form
$$
x\mathcal{D}^2_{q,x} y(x) + (\alpha+1-x)\mathcal{D}_{q,x}y(x) + [n]_qy(x) =0
$$
with solution given by the polynomial $y_{L,m_q}(x)=\sum_{\nu=0}^{n}a_{\nu}x^{\nu}$ with 
$$ a_1 = -\frac{[n]_q}{\alpha+1}a_0$$
$$a_\nu = -\frac{[n]_q}{(\alpha+1)}a_0\prod_{j=1}^{\nu-1}\left(\frac{[j]_q-[n]_q}{[j+1]_q\left([j]_q+\alpha+1\right)}\right)$$
for every $\nu\ge2$.
\end{corol}

\begin{corol}
Let $q\in\R$ with $q>0$ and $q\neq1$. For all positive integer $n$ there exists a $q$-difference equation of the form
$$\mathcal{D}_{q,x}^2y(x)-2x\mathcal{D}_{q,x}y(x)+2[n]_qy(x)=0,$$
with polynomial solution $y_{H,m_q}(x)=\sum_{\nu=0}^{n}a_{\nu}x^{\nu}$ given by
\begin{enumerate}
\item If $n$ is an even number, then $a_{2k-1}=0$ for $1\le k\le n/2$,
$$a_2=-2[n]_{q}a_0$$
and
$$    a_{2k} = -\frac{2[n]_q}{[2]_q}a_0\prod_{j=1}^{k-1}\frac{2[2j]_q-2[n]_q}{[2j+2]_q[2j+1]_q}$$
for all $k\ge 2$.
\item If $n$ is an odd number, then $a_{2k}=0$ for $0\le k\le \left\lfloor n/2\right\rfloor$, and
$$ a_{2k+1} = a_1\prod_{j=0}^{k-1} \frac{2[2j+1]_q-2[n]_q}{[2j+3]_q[2j+2]_q}$$
for $k\ge1$.
\end{enumerate}
\end{corol}

\begin{corol}
Let $q\in\R$ with $q>0$ and $q\neq1$ and $\alpha,\beta>-1$. For every integer $n\ge2$ there exists a $q$-difference equation of the form
$$(1-x^2)\mathcal{D}_{q.x}^2y(x)+[\beta-\alpha-(\alpha+\beta+2)x]\mathcal{D}_{q,x}y(x)+ \left([n]_q([n-1]_q+\alpha+\beta+2)\right)y(x)=0$$
with polynomial solution $y_{J,m_q}(x)=\sum_{\nu=0}^{n}a_{\nu}x^{\nu}$ defined by
$$a_2 = -\frac{(\beta-\alpha)a_1+[n]_q\left([n-1]_q+\alpha+\beta+2\right)a_0}{[2]_q}$$
$$a_3 = \frac{\left[\alpha+\beta+2+(\beta-\alpha)^2-[n]_q\left([n-1]_q+\alpha+\beta+2\right)\right]a_1+[n]_q\left([n-1]_q+\alpha+\beta+2\right)(\beta-\alpha)a_0}{[3]_q[2]_q}$$
and 
$$     a_{\nu+2} = \frac{[\nu]_q\left([\nu-1]_q+\alpha+\beta+2\right)-[n]_q\left([n-1]_q+\alpha+\beta+2\right)}{[\nu+2]_q[\nu+1]_q} a_\nu-\frac{\beta-\alpha}{[\nu+2]_q}a_{\nu+1}$$
for every $\nu\ge2$.
\end{corol}

\begin{corol}
Let $q\in\R$ with $q>0$ and $q\neq1$. For every integer $n\ge 1$ there exists a $q$-difference equation of the form
$$x^2D_{q,x}^2y(x)+2(x+1)D_{q,x}y(x)-[n]_q([n-1]_q+2)y(x)=0$$
with polynomial solution $y_{B,m_q}(x)=\sum_{\nu=0}^{n}a_{\nu}x^{\nu}$ defined by
$$a_{1}=\frac{[n]_q([n-1]_q+2)}{2}a_0$$

$$a_2=\frac{\left(\lambdaQ{n}{n-1}-1\right)^2-1}{4[2]_q}a_0$$
and 
$$a_{\nu}=\frac{\left(\lambdaQ{n}{n-1}-1\right)^2-1}{4[2]_q}a_0\displaystyle\prod_{j=2}^{\nu-1}\frac{\lambdaQ{n}{n-1}-\lambdaQ{j}{j-1}}{2[j+1]_q}$$
for every $\nu\ge 3$.
\end{corol}

\subsection{Numerical confluence results}

In this section, we aim to show numerically how the polynomial solutions to moment differential can approach the classical polynomials when choosing adequate sequences of moments. 

\begin{theo}
Let $(\mu_p)_{p\ge1}$ be a sequence of positive rational numbers such that $\mu_p\to1$ when $p\to\infty$. Then, $y_{\star,m_{\mu_p}}(x)$ approaches to the corresponding classical polynomial, for $\star\in\{L,H,J,B\}$, and any fixed $x\in\R$. At the same time, $y_{\star,m_{\mu_p}}$ satisfies a fractional differential equations which converges to the classical second-order differential equation satisfied by the corresponding classical polynomial.

Let $(q_p)_{p\ge1}$ be a sequence of positive real numbers with $q_p\neq 1$ for $p\ge1$, and assume that $q_p\to 1$ when $p\to\infty$. Then, $y_{\star,m_{q_p}}(x)$ approaches to the corresponding classical polynomial, for $\star\in\{L,H,J,B\}$, and every $x\in\R$. At the same time, $y_{\star,m_{\mu_p}}$ satisfies a second order $q$-difference equation which converges to the classical second-order differential equation satisfied by the corresponding classical polynomial.
\end{theo}

\begin{proof}
Most of the statements are direct from the construction of the polynomials $y_{\star,\mu_p}(x)$ and $y_{\star,q_p}(x)$. We now prove that such polynomials approach the corresponding classical ones. This can be done under more general assumptions, by considering a sequence $(m^{(p)})_{p\ge1}$ with $m^{(p)}=(m^{(p)}_n)_{n\ge0}$ being a sequence of positive real numbers such that $m^{(p)}_n\to n!$ for $p\to\infty$. Indeed, let us fix $n\ge1$ and write
$$y_{\star,m^{(p)}}(x)=\sum_{\nu=0}^{n}\frac{a_{\nu}^{\star}}{m^{(p)}_\nu}x^\nu$$
for the polynomial constructed in Proposition~\ref{prop1} when $\star=L$, resp. Proposition~\ref{prop2} when $\star=H$, resp. Proposition~\ref{prop3} when $\star=J$, resp. Proposition~\ref{prop4} when $\star=B$. The construction of the coefficients $a_{\nu}^{\star}$ in the corresponding results converge to those determined by the classical polynomials when $p\to\infty$. Taking this and the fact that $m^{(p)}_n\to p!$ whenever $p\to\infty$ yields the conclusion.
\end{proof}  

We support the previous result with numerical experiments.

\begin{figure}
	\centering
		\includegraphics[width=0.45\textwidth]{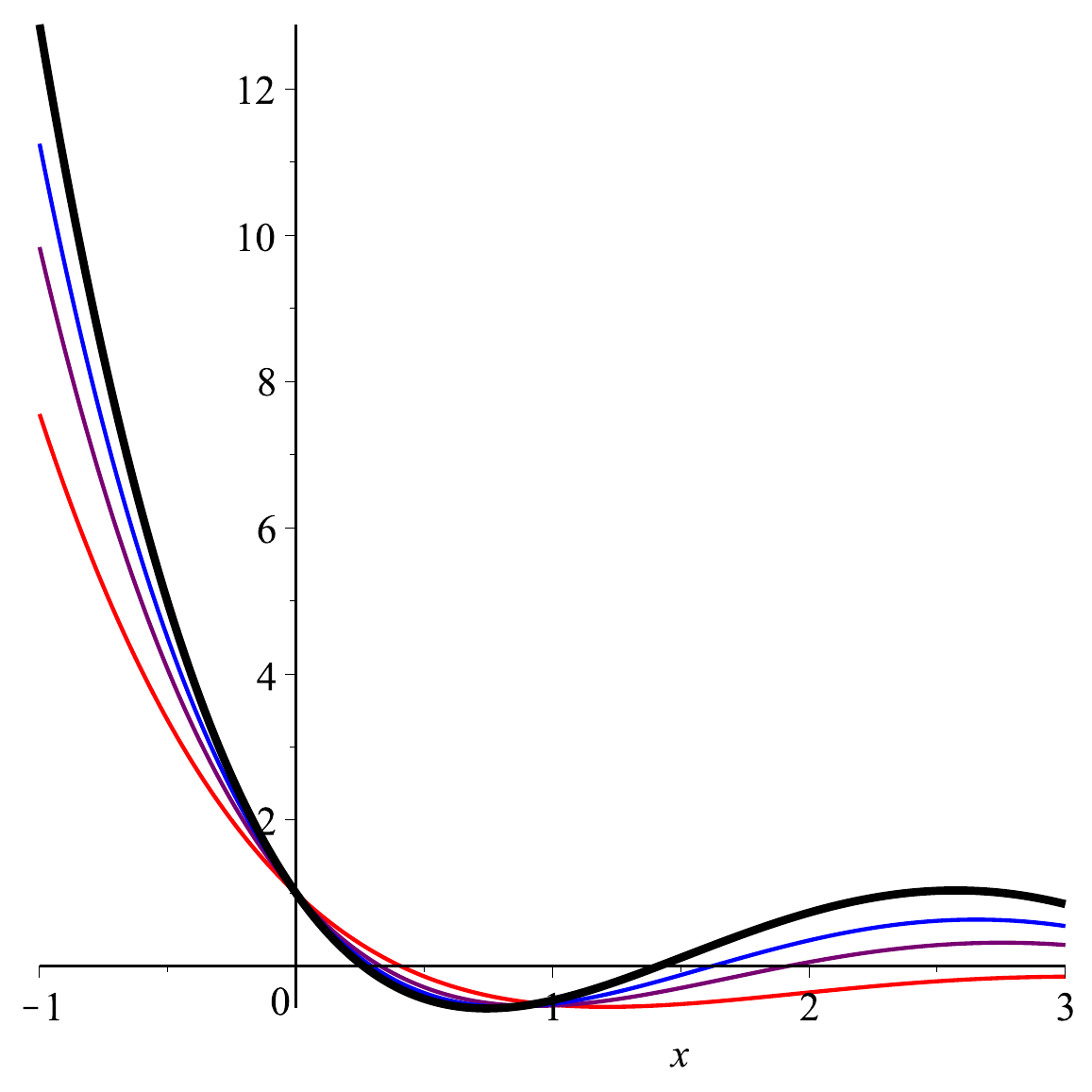}
		\includegraphics[width=0.45\textwidth]{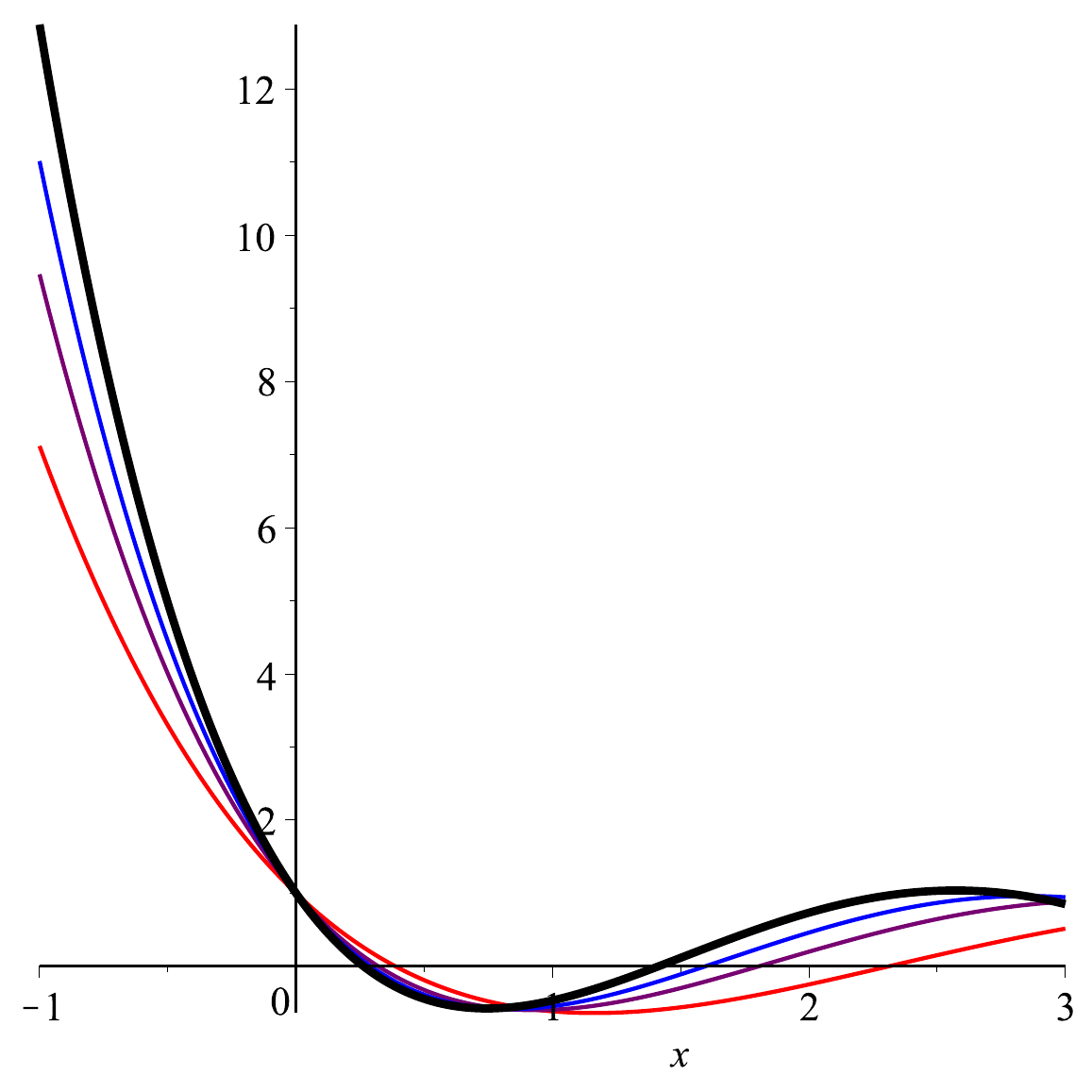}
		\caption{Approximation to Laguerre polynomial, $n=5$}
		\label{fig1}
\end{figure}

Figure~\ref{fig1} (left) displays the polynomials in $x^{\mu}$ of degree $n=5$ and $\alpha=0$ associated with classical Laguerre polynomials, which satisfy a second order fractional differential equation for the values $\mu=0.8$ (in red color), $\mu=0.9$ (purple), $\mu=0.95$ (blue). The constant term is chosen to be 1. The classical Laguerre polynomial of degree $n=5$ is drawn with a thick black curve. Figure~\ref{fig1} (right) shows the polynomials of degree $n=5$ associated with classical Laguerre polynomials, which satisfy a second order $q$-difference equation for the values $q=0.8$ (in red color), $q=0.9$ (purple), $q=0.95$ (blue). The constant term is chosen to be 1. The classical Laguerre polynomial of degree $n=5$ is drawn with a thick black curve.

\begin{figure}
	\centering
		\includegraphics[width=0.45\textwidth]{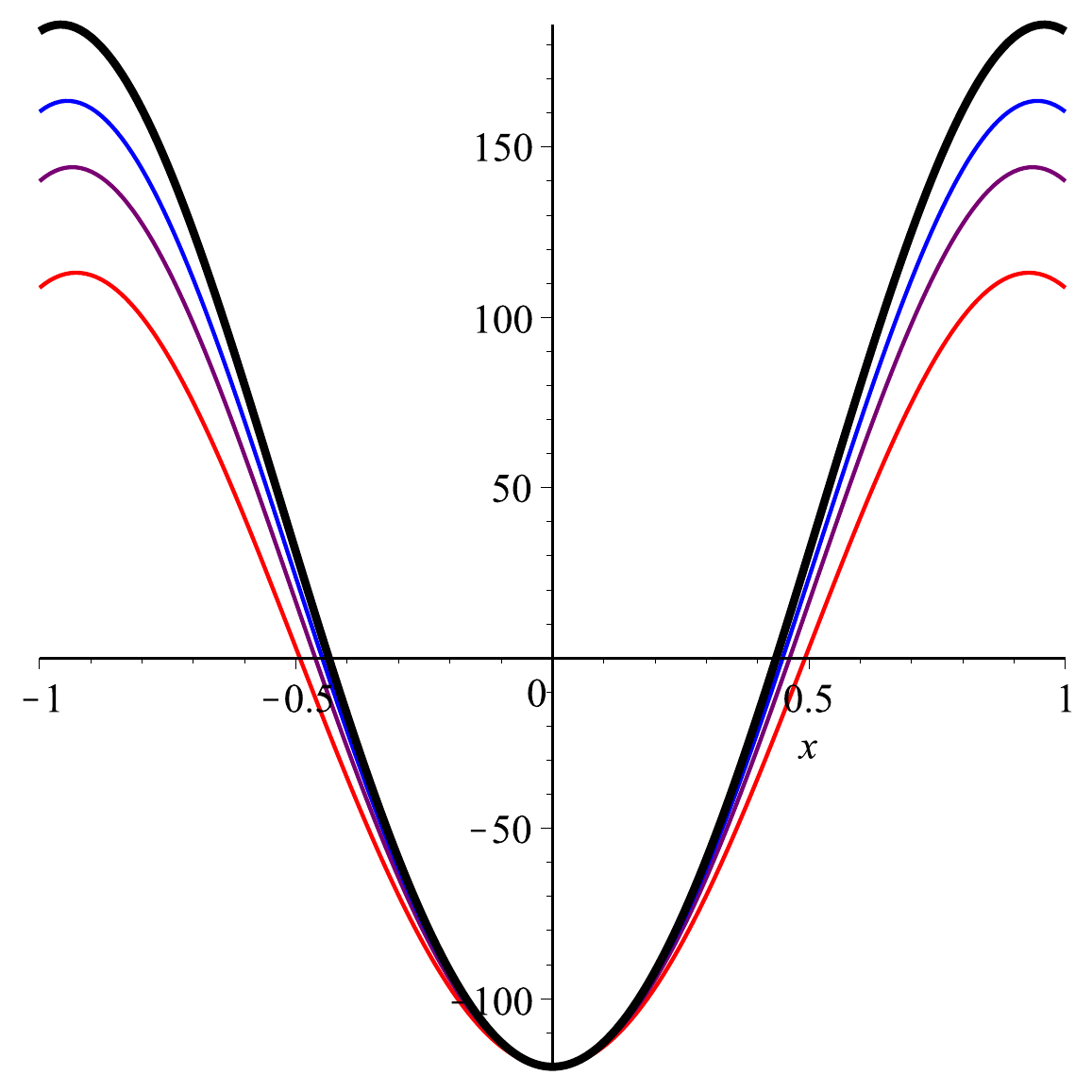}
		\includegraphics[width=0.45\textwidth]{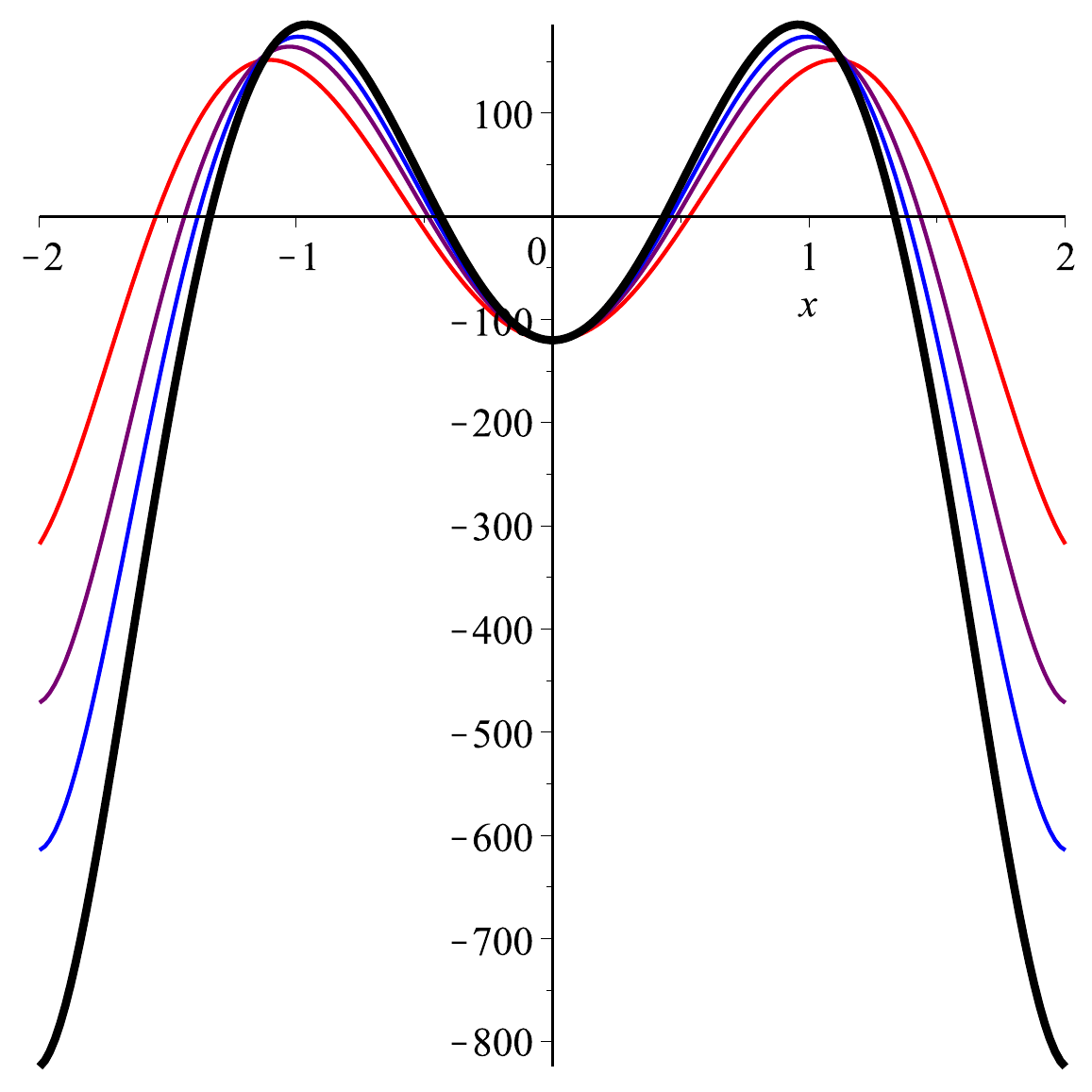}
		\caption{Approximation to Hermite polynomial, $n=6$}
		\label{fig2}
\end{figure}

In Figure~\ref{fig2} (left) we show the polynomials in $x^{\mu}$ of degree $n=6$, which coincide at $x=0$ with Hermite polynomial, and which approximate Hermite polynomial at the same time they satisfy a second order fractional differential equation for the values $\mu=0.8$ (in red color), $\mu=0.9$ (purple), $\mu=0.95$ (blue). Classical Hermite polynomial is drawn in black. Figure~\ref{fig2} (right) shows the polynomials of degree $n=6$ associated with classical Hermite polynomials, which satisfy a second order $q$-difference equation for the values $q=0.8$ (in red color), $q=0.9$ (purple), $q=0.95$ (blue). Their value at $x=0$ coincide with that of Hermite polynomial.  Classical Hermite polynomial of degree $n=6$ is drawn in black.

\begin{figure}
	\centering
		\includegraphics[width=0.45\textwidth]{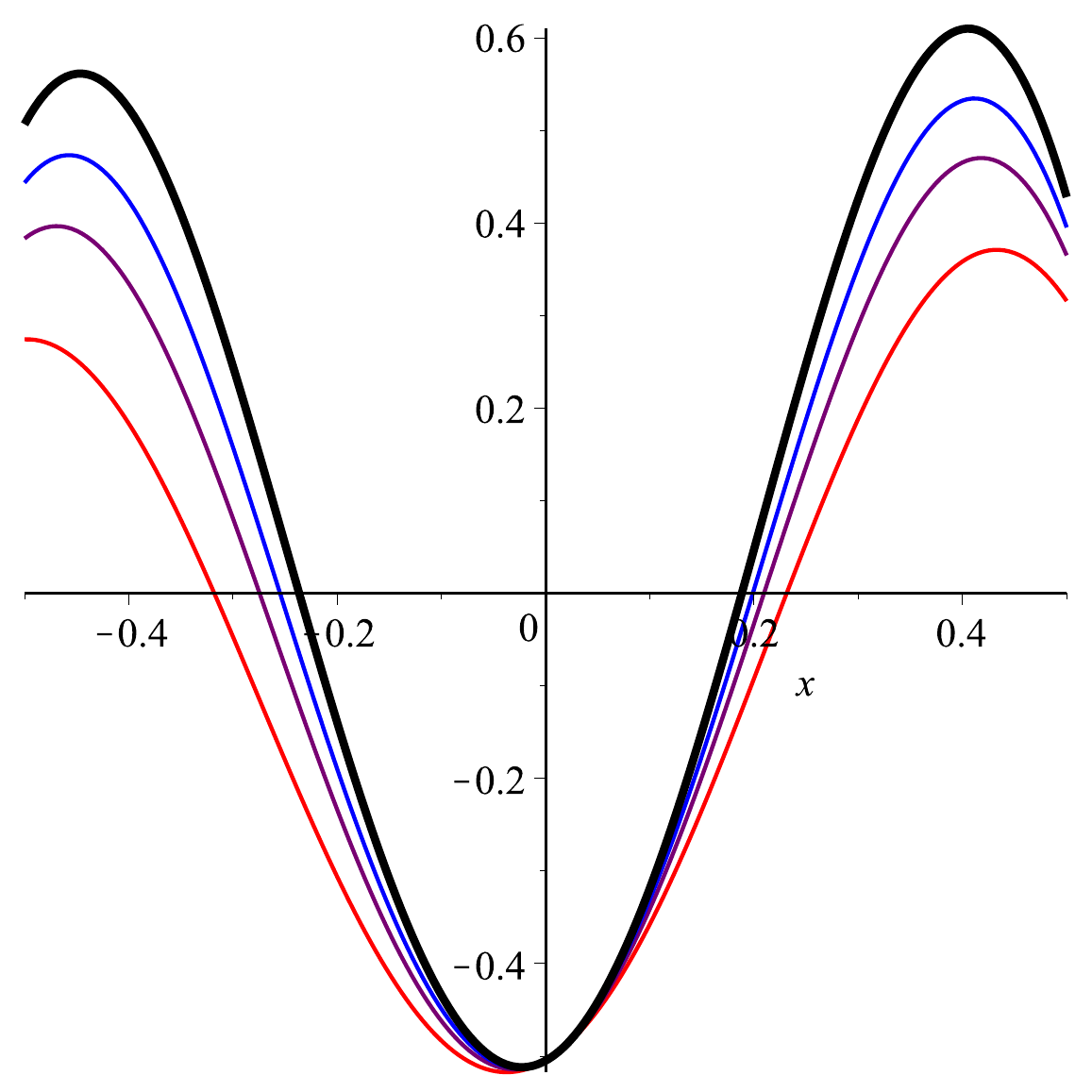}
		\includegraphics[width=0.45\textwidth]{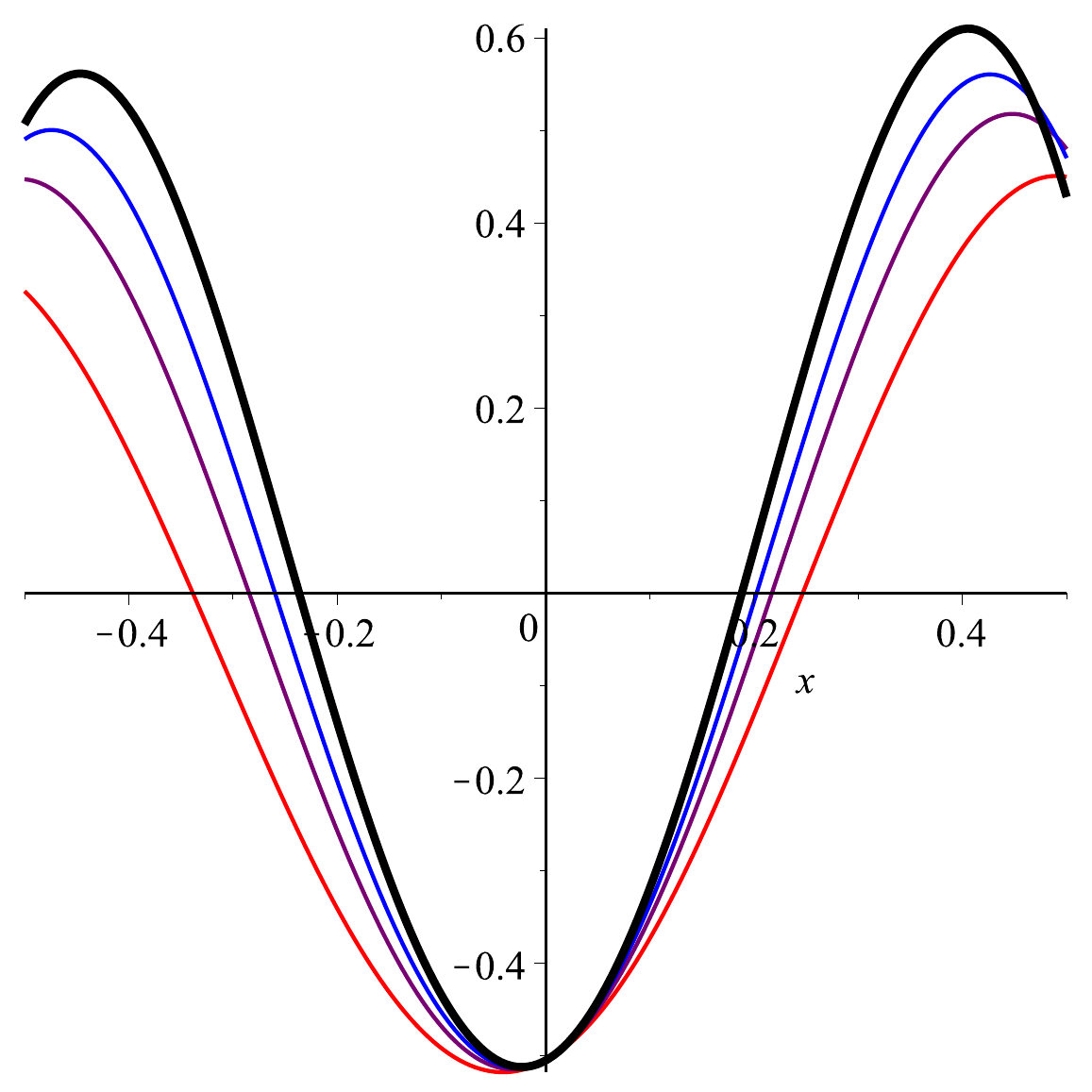}
		\caption{Approximation to Jacobi polynomial, $n=6$}
		\label{fig3}
\end{figure}

Let $\alpha=1$ and $\beta=3/4$. Figure~\ref{fig3} illustrates the approximation made by the polynomials in $x^{\mu}$ of degree $n=6$, which coincide together with their derivative with Jacobi polynomial at $x=0$. We consider the values $\mu=0.8$ (red), $\mu=0.9$ (purple), $\mu=0.95$ (blue) which approach the classical Jacobi polynomial of degree $n=6$, in black. Figure~\ref{fig2} (right) shows the approximation made with polynomials of degree $n=6$ satisfying a second order $q$-difference equation for $q=0.8$ (red), $q=0.9$ (purple), $q=0.95$ (blue), approaching classical Jacobi (black) of degree $n=6$. The polynomials are chosen so that they coincide, together with their first derivative, at $x=0$ with Jacobi polynomial of degree $n=6$.

\begin{figure}
	\centering
		\includegraphics[width=0.45\textwidth]{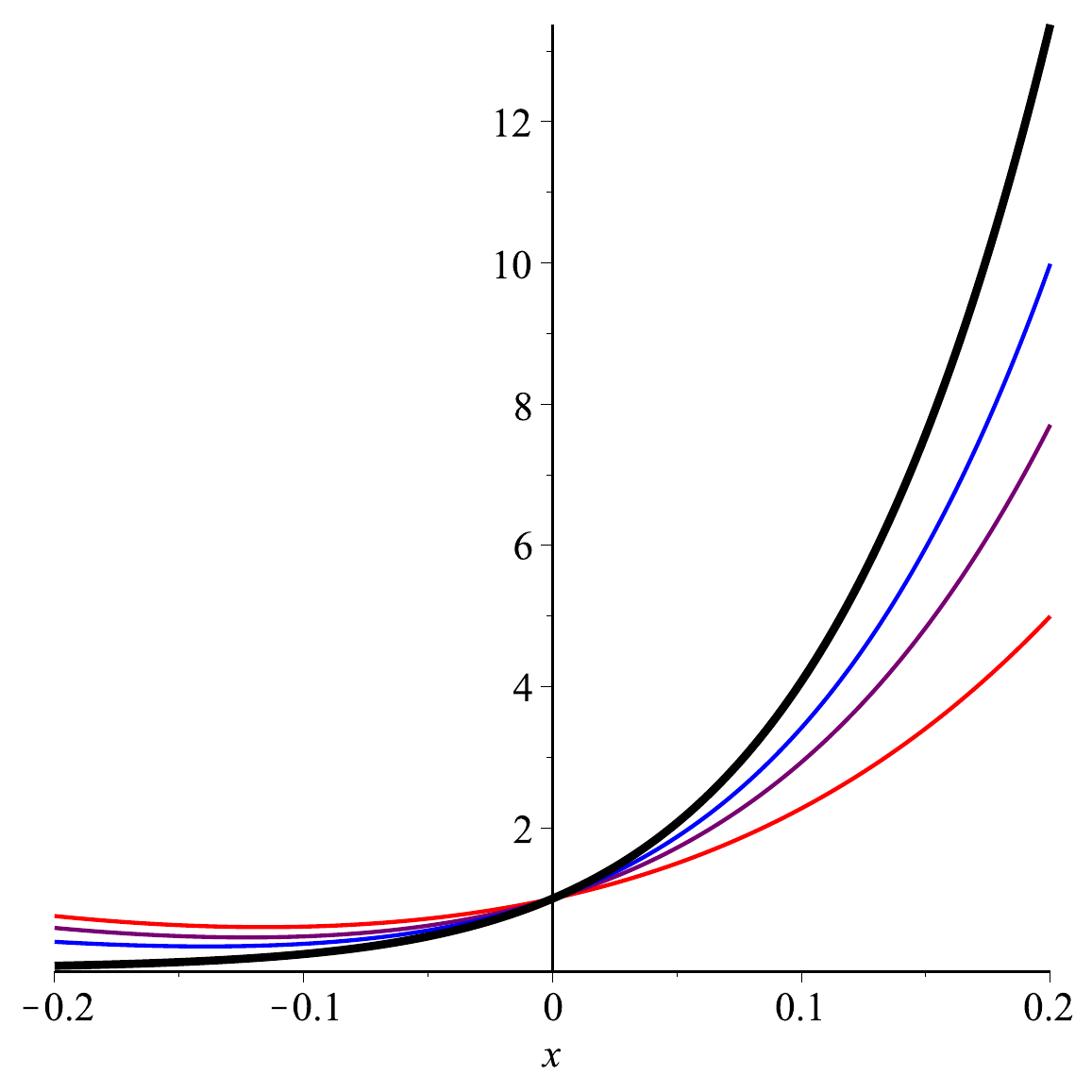}
		\includegraphics[width=0.45\textwidth]{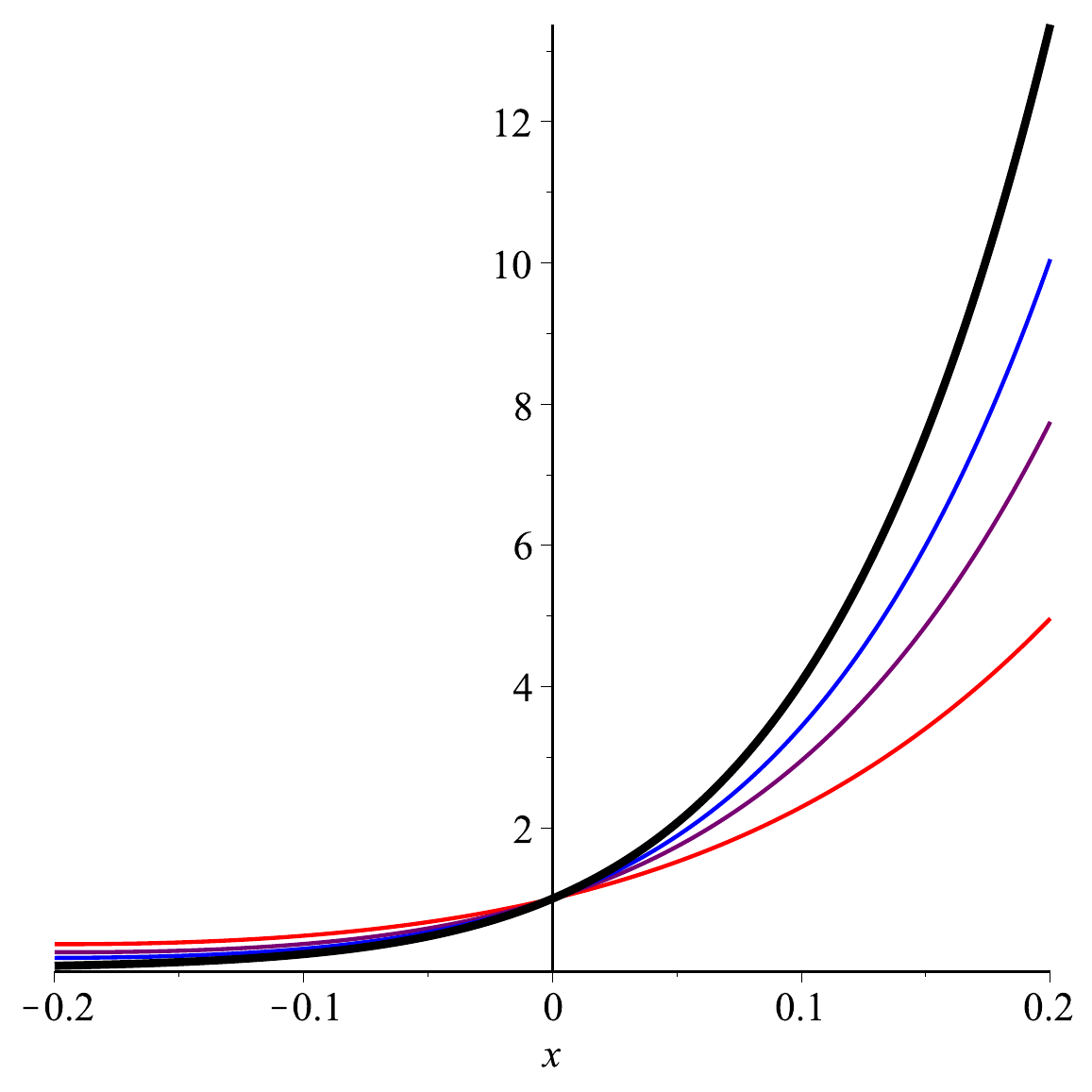}
		\caption{Approximation to Bessel polynomial, $n=5$}
		\label{fig4}
\end{figure}

In Figure~\ref{fig4} (left) we display the polynomials in $x^{\mu}$ of degree $n=5$, which coincide at $x=0$ with Bessel polynomial, and which approximate Bessel polynomial and satisfy a second order fractional differential equation for the values $\mu=0.8$ (red), $\mu=0.9$ (purple), $\mu=0.95$ (blue). Classical Bessel polynomial is drawn in black. Figure~\ref{fig4} (right) shows the polynomials of degree $n=5$ associated with classical Bessel polynomials, which satisfy a second order $q$-difference equation for the values $q=0.8$ (in red color), $q=0.9$ (purple), $q=0.95$ (blue). Their value at $x=0$ coincide with that of Bessel polynomial. We observe how they approach the classical Bessel polynomial of degree $n=5$, in black.

The confluence of the zeros of the polynomials can also be stated. As a consequence of the convergence of the coefficients of the generalized polynomials in Section~\ref{sec3} to the coefficients of the classical polynomials, one derives the convergence of their roots. More precisely, for every $n\ge1$ let $m^{(n)}=(m_p^{(n)})_{p\ge0}$ be a family of sequences of positive real numbers such that $m^{(n)}_p\to p!$ when $n\to\infty$, for every $p\ge0$. Then, for any fixed $N\ge1$ let $y$ denote the corresponding classical (Laguerre, Hermite, Jacobi or Bessel) polynomial of degree $N$. Then, there exists a sequence $(x_{n})_{n\ge1}$ of roots of the corresponding generalized polynomial of degree $N$, say $y_{m^{(n)}}$, i.e. $y_{m^{(n)}}(x_n)=0$ for every $n\ge 1$, such that $x_{n}\to x_0$, where $x_0$ is a root of $y$. See Proposition 5.2.1~\cite{artin}, for a proof of the previous statement.

Moreover, the explicit nature of the coefficients of the generalized polynomials constructed in Section~\ref{sec3} allows to give upper bounds on the distance of the roots of the generalized polynomials and the classical ones in virtue of the following classical result.

\begin{theo}[Theorem 1,~\cite{bek}]
Let 
$$f(z)=z^n+a_1z^{n-1}+\ldots+a_n=\prod_{i=1}^{n}(z-\alpha_i),$$
$$g(z)=z^n+b_1z^{n-1}+\ldots+b_n=\prod_{i=1}^{n}(z-\beta_i).$$
Then, there exists a permutation of the roots of $f$ and $g$ such that after such permutation one has 
$$\max_{i}|\alpha_i-\beta_i|\le2^{2-1/n}\left(\sum_{k=1}^n|a_k-b_k|\gamma^{n-k}\right)^{1/n},$$
where $\gamma=2\max_{1\le k\le n}\left(|a_k|^{1/k},|b_k|^{1/k}\right)$. 
\end{theo}

This result refines the previous, not only stating convergence but also a convergence rate of the roots for each of the families of polynomials. 

\begin{example}
We consider the generalized Laguerre polynomials with $\alpha=0$, associated to the sequence $m_q=([p]_q!)_{p\ge0}$. Taking into account the coefficients of $y_{L,m_q}$ of degree $n$ obtained in Section~\ref{sec42}, we have that the difference of corresponding roots of $y_{L,m_q}$ and the classical Laguerre polynomial of degree $n$ is upper bounded by
$$2^{2-1/n}\left[\sum_{k=1}^{n}\left|\prod_{j=n-k}^{n-1}\frac{[j+1]_q([j]_q+1)}{[j]_q-[n]_q}-\frac{n!^2}{(n-k)!^2(-1)^{k-1}k!}\right|\gamma^{k}\right]^{1/n},$$
with 
$$\gamma=2\max_{1\le k\le n}\left\{\left(\prod_{j=n-k}^{n-1}\frac{[j+1]_q([j]_q+1)}{[j]_q-[n]_q}\right)^{1/k},\left(\frac{n!^2}{(n-k)!^2(-1)^{k-1}k!}\right)^{1/k}\right\}.$$
\end{example}
\section{Conclusions}

For each of the families of classical orthogonal polynomials: Laguerre, Hermite, Jacobi and Bessel, we construct for every non-negative integer $n$
\begin{enumerate}
\item a second-order moment differential equation,
\item a polynomial of degree $n$,
\end{enumerate}
generalizing the classical equations and polynomials. Two particular manifestations of the sequence of moments, of great importance in applications are considered. In a first application, we are constructing for every $n$ a second-order fractional differential equation and a polynomial-like solution of degree $n$. The adequate modification of the fractional order allows to provide a confluence of the fractional differential equation to the classical second-order differential equation satisfied by the corresponding classical polynomial, whilst the polynomial-like solution tends to the classical polynomial. In a second application, we are constructing for every $n$ a second order $q$-difference equation and a polynomial solution of degree $n$, for every $q>0$. If one considers $q\to 1$, then the $q$-difference equation approaches to the classical second-order differential equation satisfied by the corresponding classical polynomial, and the polynomial solution of the $q$-difference equation tends to the classical orthogonal polynomial. 

Many questions are to be answered in a future direction, such as the existence (or not) of a positive measure for which the polynomials obtained are orthogonal. In an affirmative case, is the sequence of moments considered related to the moments of such measure? In this concern, this should be considered as a seminal work to be continued in a future research to answer these and other questions.   
\section*{Acknowledgements}

The work of E. J. Huertas and A. Lastra has been supported by Direcci\'{o}n General de Investigaci\'{o}n e Innovaci\'{o}n, Consejer\'{i}a de Educaci\'{o}n e Investigaci\'{o}n of the Comunidad de Madrid (Spain) and Universidad de Alcal\'{a}, under grant CM/JIN/2021-014, \textit{Proyectos de I+D para J\'{o}venes Investigadores de la Universidad de Alcal\'{a} 2021}, and the Ministerio de Ciencia e Innovaci\'{o}n-Agencia Estatal de Investigaci\'{o}n MCIN/AEI/10.13039/501100011033 and the European Union ``NextGenerationEU''/PRTR, under grant TED2021-129813A-I00.\\

This research was conducted while E. J. Huertas was visiting the ICMAT (Instituto de Ciencias Matemáticas), from jan-2023 to jan-2024 under the Program \textit{Ayudas de Recualificación del Sistema Universitario Español para 2021-2023 (Convocatoria 2022) - R.D. 289/2021 de 20 de abril (BOE de 4 de junio de 2021)}. This author wish to thank the ICMAT, Universidad de Alcal\'{a}, and the Plan de Recuperaci\'{o}n, Transformaci\'{o}n y Resiliencia (NextGenerationEU) of the Spanish Government for their support.\\

The work of A. Lastra is also partially supported by the project PID2019-105621GB-I00 of Ministerio de Ciencia e Innovación, Spain.\\

The work of V. Soto-Larrosa has been supported by Consejer\'ia de Econom\'ia, Hacienda y Empleo of the Comunidad de Madrid through ``Programa Investigo'', funded by the European Union ``NextGenerationEU''.\\

The authors are members of the research group AnFAO (Cod.: CT-CE2023/876) of Universidad de Alcalá.

\end{document}